\documentclass[12pt]{article}
\usepackage{amssymb}
\usepackage{latexsym}
\usepackage{exscale}

\usepackage{pgf,tikz}

\usetikzlibrary{trees}

\usepackage{amssymb, amsmath}
\usepackage{enumerate}
\usepackage{enumerate}
\usepackage[top=1in,left=1in,right=1in,bottom=1in]{geometry}

\usepackage{xcolor}
\usepackage{soul}

\usepackage{cancel}

\def \beq{\begin{equation}}
\def \eeq{\end{equation}}
\def \lab{\label}
\renewcommand{\rq}[1]{(\ref{#1})}
\newtheorem{lemma}{Lemma}
\newtheorem{prop}{Proposition}
\newtheorem{thm}{Theorem}
\newtheorem{cor}{Corollary}

\newcommand{\bR}{{ \mathbb R  }}
\newcommand{\bC}{\Bbb C}
\newcommand{\bZ}{\Bbb Z}
\newcommand{\bN}{\Bbb N}

\newcommand{\bQ}{\Bbb Q}

\newcommand{\la}{\mbox{$\lambda$}}

\newcommand{\pa }{\partial }

\newcommand{\ep}{\epsilon}

\newcommand{\al}{\alpha }
\newcommand{\be}{\beta }

\newcommand{\f}{\varphi }

\newcommand{\mH}{{\mathcal H}}
\newcommand{\iy}{{\infty}}
\newcommand{\om}{{\omega}}
\newcommand{\Om}{{\Omega}}

\newcommand{\ka}{{\kappa}}

\def\<{\langle} \def\>{\rangle}

\begin{document}

\title{Shape, Velocity, and Exact Controllability for the Wave Equation on a Graph with Cycle}

 \maketitle

\noindent{\bf Sergei Avdonin}, Department of Mathematics and Statistics, University of Alaska Fairbanks, Fairbanks,
AK 99775, U.S.A.; s.avdonin@alaska.edu;  \\ and Moscow Center for Fundamental and Applied Mathematics,  Moscow 119991, Russia

\

\noindent{\bf Julian Edward}, Department of Mathematics and Statistics, Florida International University, Miami,
FL 33199, U.S.A.; edwardj@fiu.edu; 

\

\noindent{\bf Yuanyuan Zhao}, Department of Mathematics and Statistics, University of Alaska Fairbanks, Fairbanks,
AK 99775, U.S.A.; yuanyuanzhao17@gmail.com; 

\

{\bf Abstract.} Exact controllability is proven on a graph with cycle. The controls can be a mix of controls applied at the boundary and interior vertices. The method of proof first uses a dynamical argument to prove shape controllability and velocity controllability, thereby
solving their associated moment problems. This enables one to solve the moment problem associated to exact controllability.
In the case of a single control, either boundary or interior, it is shown that exact controllability fails.

\vskip1cm  Dedicated to memory of Sergey Naboko, a brilliant mathematician and a long time friend of the first author of this paper.

\vskip1cm 

\begin{section}{Introduction}

Controllability properties of the wave equation is a central topic of the control theory of partial differential equations. A large number of { papers describe } many powerful methods, which prove controllability of the wave equation in various spatial domains under the action of various types of controls (see, e.g. \cite{Rus}, \cite{Li}, \cite{LT}, \cite{Zu} and references therein). In this paper we describe an approach that is based on the 
 relationship between exact controllability, on one hand, and shape and velocity controllability on the other hand. This relationship was used in \cite{ABI} for a vector wave equation, in \cite{AE}, \cite{AE1} for a string with attached point masses, and in \cite{AZ} for the wave equation on a metric tree graph. Here we consider a control problem for the wave equation on a graph with cycle --- a ring with attached edge which is called the lasso graph.

Control problems for the wave equation on graphs have important applications in science
and engineering and were studied in many papers (see the monographs \cite{AI}, \cite{DZ},   \cite{LLS}; the
surveys \cite{A2}, \cite{Z}; and references therein). They also have deep connection with inverse problems on graphs, see, e.g. \cite{Be}, \cite{AK}, \cite{AAAE}, \cite{AE3}, \cite{Ku2}. 
In this paper we consider exact controllability for the wave equation of the form
$$u_{tt} - u_{xx} + q(x)u = 0.$$
There is a growing body of work in the case where  the graph is a tree, i.e. a graph without cycles, and the controls are assumed to act on the boundary. Typically, the so-called Kirchhoff-Neumann (KN) conditions
are assumed at all interior vertices. This problem was studied, e.g. in \cite{BV}, \cite{DZ}, \cite{LLS},    \cite{AZ} (in
those papers the problem was stated in slightly different forms). It was proved that the
system is exactly controllable if the control functions act at all or at all but one of the
boundary vertices.

In the case of graphs with cycles, we are unaware of any positive results concerning exact controllability. Results on lack of boundary exact controllability for graphs with cycles have been proved in  \cite[Theorem VII.5.1]{AI}, see also \cite[Remark 6.11]{DZ}. In these chapters the authors consider the case $q=0,$ but the results can be easily extended to general $q.$  Lack of boundary controllability for some graphs with cycles have been also proved in \cite[Sections II.5, V.2]{LLS}. We also mention another concept of controllability, the so-called nodal profile controllability, where at a given vertex the state and the velocity profile is exactly matched by boundary controls. This concept is weaker than the exact controllability discussed in this article, but admits cycles (see \cite{LRW}, \cite{ZLL} for further information and references).
		To reach the exact controllability of systems on graphs with cycles we need to use not only boundary but also interior controls, as was proposed in \cite{A5}.

In this paper, we consider two different control problems,  both including interior controls, for the graph, denoted $\Omega$, consisting of a ring with attached edge, see Figure \ref{lasso}. 
We note that the inverse problem for a magnetic Schr\"odinger equation on such a graph was studied in \cite{Ku1}, \cite{Ku}.

\begin{figure}[h]
\begin{center}
\begin{tikzpicture}[line cap=round,line join=round,x=0.5cm,y=0.5cm]
\draw[color=black,  ](22,5) circle (4);
\draw    (12,5) -- (18,5) ;
		\draw (18,5) node [anchor=north west][inner sep=0.75pt]   [align=left] {$x=0$};
		
		\draw (12,5) node [anchor=north west][inner sep=0.75pt]   [align=left] {$x=l$};
			\draw (26.25,5) node [anchor=north west][inner sep=0.75pt]   [align=left] {$x=a$};
		\draw (22.25,9) node [anchor=north west][inner sep=0.75pt]   [align=left] {$e_2$};
		\draw (22.25,1) node [anchor=north west][inner sep=0.75pt]   [align=left] {$e_3$};
		\draw (15,4.75) node [anchor=north west][inner sep=0.75pt]   [align=left] {$e_1$};
	\draw [shift={(18,5)}, rotate = 0] [color={rgb, 255:red, 0; green, 0; blue, 0 }  ][fill={rgb, 255:red, 0; green, 0; blue, 0 }  ][line width=0.75]      (0, 0) circle [x radius= .1, y radius= .1]   ;
	\draw [shift={(12,5)}, rotate = 0] [color={rgb, 255:red, 0; green, 0; blue, 0 }  ][fill={rgb, 255:red, 0; green, 0; blue, 0 }  ][line width=0.75]      (0, 0) circle [x radius= .1, y radius= .1]   ;
	\draw [shift={(26,5)}, rotate = 0] [color={rgb, 255:red, 0; green, 0; blue, 0 }  ][fill={rgb, 255:red, 0; green, 0; blue, 0 }  ][line width=0.75]      (0, 0) circle [x radius= .1, y radius= .1]   ;
\end{tikzpicture}
\caption{The lasso graph parametrized as three edge star. } \label{lasso}
\end{center}
\end{figure}
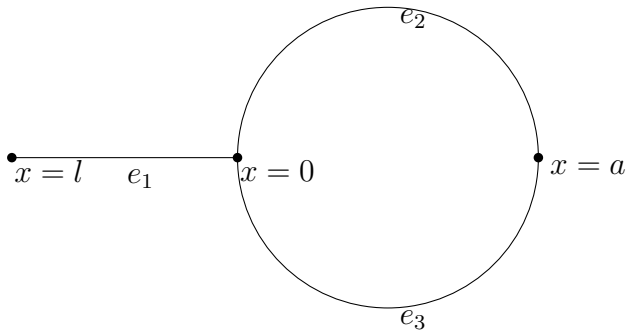

First we prove the shape and velocity controllability using the dynamical method ---  we reduce these problems to the  Volterra integral equations of the second kind. Then we prove exact controllability using the spectral approach --- the method of moments and properties of exponential families.
 This approach was used in \cite{AZ} for tree graphs, with Dirichlet boundary controls. 
In the present paper  we use both Dirichlet and Neumann type controls. This makes the spectral part of the proof more complicated.

To state our control problems, we first establish some notation. We identify the boundary edge, $e_1$,  with the interval $(0,l)$ with $x=l$ corresponding to the boundary. The cycle we will identify in two different ways. At times, it will be 
convenient to consider one more (artificial) interior vertex in the middle of the ring where we put the homogeneous Kirchoff--Neumann conditions,  and each of the  two interior edges, $e_2,e_3,$ with the interval $(0,a)$, see Figure 1. At other times, it will be convenient to denote by $e_r$ the edge starting and finishing at $x=0.$

We now state the two initial boundary value problems (IBVP) associated to this graph and two results about exact controllability of the corresponding wave equations.
In what follows, the restriction of a function $u$ to edge $e_j$ will be denoted $u_j$ and
 $\partial u_j$ will denote the derivative of $u$ at a vertex, along the incident edge $e_j$ in the direction outward from the vertex. The potential $q$ is a real valued function such that $q_j:=q|_{e_j} \in L^1(e_j).$

{\bf Problem 1.}  We consider the  IBVP with one control acting at the boundary vertex and another control acting at the interior vertex :
\beq \lab{eqq1}
u_{tt} - u_{xx} + q(x) u =0,\ t>0,\ x\in \cup_{j=1}^3 \, e_j , 
\eeq
\beq \label{art1}
u_2(a,t)=u_3(a,t); \ \pa u_2(a,t)+\pa u_3(a,t)=0, \ t>0,
\eeq 
\beq \label{ic1}
u(x,0)=u_t(x,0)=0, \ x\in \Omega , 
\eeq 
\beq \lab{bcc1} 
\pa  u_1(l,t)=f_1(t),\ t>0,
 \eeq
 \beq \lab{inc1}
u_2(0,t)-u_1(0,t)=f_2(t), \ u_3(0,t)-u_1(0,t)=0, \ \sum_{j=1}^3 \partial u_j(0,t)=0, \ t>0.  
\eeq
In what follows, we will refer to $f_1$ as a boundary control and refer to $f_2$ as an interior control.
	
	For example, in electric networks the conditions \rq{inc1} mean that a source of energy is established at the interior vertex which provides the difference of voltage between $e_2$ and two other edges, and the Kirchhoff law is satisfied for currents.

To define our relevant function spaces, we denote the cycle edge by $e_r=(0,2a)$. We denote the space $L^2(\Omega) $ by $ \mH $. We equip the space 
$$ 
 \mH^1=\{ \phi \in \mH: \ \phi_1 \in H^1(e_1), \ \phi_r \in H^1(e_r)\}
$$ 
with norm 
 $\| \phi\|^2_{\mH^1}=\| \phi '\|_{L^2(e_1)}^2+\| \phi '\|_{L^2(e_r)}^2+\| \phi \|_{\mH}^2.$ 
We then define 
$$
{\mH}^1_0=\{ \phi \in \mH^1: \phi_1(0)=\phi_r(0)=\phi_r(2a)\},\ 
{\mH}^1_1=\{ \phi \in \mH^1: \phi_1(0)=\phi_r(2a)\}.
$$
The space $\mH_1^1$ will be needed to discuss well-posedness of Problem 1.

The following theorem concerning controllability of the system \rq{eqq1}--\rq{inc1} presents the main result of the paper. 

\begin{thm}\label{thm1}
Let $T_*=a+l$ , and let $T\geq 2T_*$.
Let $(\phi_1,\phi_2) \in \mH^1_0 \times \mH$.
 There exist controls ${\bf f}:= (f_1,f_2)\in  L^2(0,T) \times H^{1}_0(0,T)$ such that  
 the solution $u^{{\bf f}}$ to the IBVP \rq{eqq1}--\rq{inc1} satisfies the equalities
\beq
u^{{\bf f}}(x,T)=\phi_1(x),\ u_t^{{\bf f}}(\cdot,T)=\phi_2(\cdot),\label{ec}
\eeq
and there exists a constant $C$ that depends only on $q,a,l$
such that 
\beq 
\| f_1\|_{L^2(0,T)}+\| f_2\|_{H^1(0,T)}\leq C(
\| \phi_1\|_{\mH^1}+\| \phi_2\|_{\mH}).\label{ce1}
\eeq
\end{thm}

The existence of ${\bf f}$ satisfying \rq{ec}, \rq{ce1} is called exact controllability. The system would be ``approximately controllable"   (a weaker conclusion)  if for any
$\ep >0$, there exists ${\bf f}$ satisfying \rq{ce1} and (instead of \rq{ec}) 
$$\| u^{{\bf f}}(\cdot,T)-\phi_1\|_{\mH^1}+\| u_t^{{\bf f}}(\cdot,T)-\phi_2\|_{\mH}<\ep .
$$

\ 

{\bf Problem 2.} Now we consider the same graph with three controls at the interior vertex. Thus in addition to \rq{eqq1}, \rq{art1}, and \rq{ic1}, we have 
\beq \lab{bc2} 
\pa  u_1(l,t)=0, \ t>0,
 \eeq
\beq \lab{3c1}
u_3(0,t)-u_1(0,t)=f_3(t), \ u_2(0,t)-u_1(0,t)=f_2(t), \ \sum_{j=1}^3 \partial u_j(0,t)=f_1(t),  
\ t>0.
\eeq
Let $T^*=\max (a,l).$

\begin{thm}\label{thm2}
Let $T> 2T^*$.
Let $(\phi_1,\phi_2) \in \mH^1_0 \times \mH$.
 There exist controls 
 $${\bf f}:= (f_1,f_2,f_3)\in L^2(0,T)\times H^1_0(0,T)\times H^1_0(0,T) $$
 such that  the solution $u^{{\bf f}}$ to the IBVP \rq{eqq1}--\rq{ic1}, \rq{bc2}, \rq{3c1} satisfies the equalities
$$u^{{\bf f}}(x,T)=\phi_1(x),\ u_t^{{\bf f}}(\cdot,T)=\phi_2(\cdot).$$
Furthermore, there exists a constant $C$ that depends only on $q,a,l, T$
such that 
\beq 
\| f_1\|_{L^2(0,2T)}+\| f_2\|_{H^1(0,2T)} +\| f_3\|_{H^1(0,2T)}\leq C(
\| \phi_1\|_{\mH^1}+\| \phi_2\|_{\mH}).\label{ce2}
\eeq
\end{thm}

We view Problems 1 and 2 as toy models, that can be applied to more general graphs with cycles. This will be done in future work. 

A key ingredient in the proofs of Theorems 1,2 are shape and velocity control results that we state now just for Problem 1. 
\begin{thm}\label{thm3}

\

a) (Shape control) Let $T\geq T_*$. Let $\phi_1 \in\mH^1_0$.
 There exist controls ${\bf f}:= (f_1,f_2)\in L^2(0,T)\times H_0^1(0,T)$ such that 
$$u^{{\bf f}}(x,T)=\phi_1(x),$$
and there exists a constant $C$ that depends only on $q,a,l$
such that  the solution $u^{{\bf f}}$ to the IBVP \rq{eqq1}--\rq{inc1} satisfies the equalities
\beq 
\| f_1\|_{L^2(0,T)}+\| f_2\|_{H^1(0,T)}\leq C
\| \phi_1\|_{\mH^1}.\label{ce3}
\eeq
\

b) (Velocity control). Let $T\geq T_*$. Let $\phi_2 \in \mH$.
 There exist controls ${\bf f}:= (f_1,f_2)\in L^2(0,T)\times H^1(0,T)$ such that 
$$u_t^{{\bf f}}(x,T)=\phi_2(x),$$
and there exists a constant $C$ that depends only on $q,a,l$
such that  the solution $u^{{\bf f}}$ to the IBVP \rq{eqq1}--\rq{inc1} satisfies the equalities
\beq 
\| f_1\|_{L^2(0,T)}+\| f_2\|_{H^1(0,T)}\leq C
\| \phi_2\|_{\mH}.\label{ce4}
\eeq

\end{thm}
We remark in passing that the notions of shape controllability and velocity controllability  are of independent interest because of their connection with solvability of inverse problems, \cite{A5}.

To prove Theorems \ref{thm1},\ref{thm2}, we show that solvability of the moment problems associated to shape and velocity control implies the solvability of the moment problem associated to exact control in twice the time.

This paper is organized as follows. In Section \ref{rep} we give  representations of the solution of a wave equation with boundary control; these representations are
then used to prove the shape and velocity control results in 
Section \ref{s3}. Theorems \ref{thm1} and \ref{thm2} are then proved in Section \ref{s4}. { Then in Section \ref{s5}  we discuss lack of approximate and exact controllability
of the system with a smaller number of controls.}

\end{section}

\begin{section}{Solution representations}\label{rep}
\begin{subsection}{Solution representation on half-line}

As a first step to discuss  controllability for Problems 1 and 2, it is useful to
get a representation of the solution in the case of a single control on the half line.
 For functions of only one variable,
 $f^{(j)}$ denotes the $j$th derivative. When convenient we will denote
 $f^{(1)}=f'$, resp. $f^{(2)}=f''$. For partial derivatives of $u$ we use the notations
  $u_{xx}$ or $\frac{\partial^2u}{\partial x^2}$.

  Define $H^{j}(a,b)$ to be the set of functions in $L^2[a,b]$ whose
 weak derivatives up to order $j$ are in $L^2(a,b)$. The corresponding norms will be denoted $||*||_{H^j(a,b)}$.

As a first step to solving  Problems 1 and 2,
we will find a useful representation of the solution to the following system that models a semi-infinite vibrating string, with control at $x=0$.
\begin{equation}
u_{tt}-u_{xx}+q(x)u=0, \ t\in (0,\infty ),\ x\in (0,\infty ) ,\label{pdeR}
\end{equation}
\begin{equation}
 \beta_1 u(0,t)+\beta_2u_x(0,t)=g(t),\ t>0,\ \beta_1^2+\beta_2^2=1\label{bcR}
 \end{equation}
\begin{equation}
 u(x,0)=u_t(x,0)=0,\ x>0 .\label{icR}
 \end{equation}
In this paper we will consider only Neumann  boundary 
conditions ($\beta_2=1$), or Dirichlet boundary conditions ($\beta_1=1$).

 As a preliminary step,
 consider the Goursat-like problem in the domain $D =\{ (x,t)\ |\ 0<x<t<\infty \}$:

 \ \vspace{-10pt}
\begin{eqnarray}
k_{tt}(x,t)-k_{xx}(x,t)+q(x)k(x,t) & = & 0,\ (x,t)\in D,\label{goum}\\
\beta_1 k(0,t) + \beta_2 k_x(0,t)& = & 0,\ t>0,\label{mbcm}\\
 k(x,x)& = & -\frac{1}{2}\int_{0}^{x} {q}(\eta  )d\eta ,\ x\geq 0. 
\label{g2m}
\end{eqnarray}

\begin{prop} \label{goursatm}
\

A)   Let ${q}\in C[0,\infty )$.
 System \eqref{goum}-\eqref{g2m} has a unique generalized
  solution, denoted $k_{}( x,t)$, such that $k_{}(\cdot,\cdot)\in C^{1}(\overline {D} )$ and the boundary conditions hold in a classical sense.

B)  Let ${q}\in L^1_{loc}[0,\infty )$.
 System \eqref{goum}-\eqref{g2m} has a unique  generalized 
  solution, denoted $k_{}( x,t)$, and 
\beq \lab{reg1}
  k_x(\cdot,s),\ k_s(\cdot , s),\ k_x(x,\cdot ), \ k_x(x,\cdot )\in L^1_{loc}(0,\infty).
  \eeq
  The partial derivatives appearing in \rq{reg1} depend continuously in $L^1_{loc}(0,\infty)$ on the parameters $s,x,$
  and \eqref{goum} and the boundary conditions hold almost everywhere.
\end{prop}
The solution to the standard Goursat problem in the case $\beta_2=0$ is well known (see, e.g. \cite{AMi}), 
and in the general case the reader is referred to \cite{AE1}.

The following now holds by direct calculation.
\begin{prop}\label{wavesolm}
Let $T>0$. Let $k$ be as in Proposition \ref{goursatm} with ${q}\in L^1_{loc}(0,\infty)$.

In all statements below, assume $g(t)=0$ for $t<0$.

A) (Neumann and mixed control, $\beta_2 \neq 0$).
 Let
\begin{equation}
f(t)=- \frac{1}{\be_2} \int_{s=0}^tg(s) \, \exp\left\{\frac{\be_1}{\be_2} (t-s)\right\}\, ds.\label{fg}
\end{equation}

 a) Suppose $g\in C^1[0,\iy)$ with $g(0)=0$. Then
System \eqref{pdeR},\eqref{bcR},\eqref{icR} has unique solution
$u^g(x,t)$, with
\begin{equation}
u^g(x,t)=
f(t-x)+\int_{s=x}^tk(x,s)f(t-s)ds, \label{wf1}
\end{equation}
$u^g\in H^{2}((0,\infty)\times
(0,T))$, \eqref{pdeR} is satisfied almost everywhere, and the
boundary and initial conditions are satisfied in a classical sense.

b)  For $g\in L^2(0,\iy)$, the function $u^g(x,t)$ defined above
gives a  solution to
\eqref{pdeR} in the distribution sense,  \eqref{bcR} holds almost everywhere in $t$, and \eqref{icR}  holds for all $x$.
Furthermore,   $u^g\in
C([0,T];H^1(0,\infty ))$.

B) (Dirichlet control, $\beta_2=0$).

 a) Suppose $g\in C^2[0,\iy)$ with $g(0)=g'(0)=0$. Then
System \eqref{pdeR},\eqref{bcR},\eqref{icR} has unique solution
 $u^g(x,t)$, with
\begin{eqnarray}
u^g(x,t) & = &
g(t-x)+\int_{s=x}^tk(x,s)\,g(t-s)\,ds,\label{wf1d}
\end{eqnarray}
  $u^g\in H^{2}((0,\infty)\times
  (0,T))$ for each $T>0$, the equation  \eqref{pdeR} is satisfied almost everywhere, and the
  boundary and initial conditions are satisfied in a classical sense.

 b) 
 For $g\in L^2(0,T)$, the function $u^g(x,t)$ defined above
 gives the unique  solution to the generalized boundary value problem,  and $u^g\in C(0,T;W^T)$, where 
 $W^T=\{ u\in L^2_{loc}(0,\infty): support(u)\subset [0,T]\}.$

\end{prop}

{\bf Remark 1.}
Since $g(t)=0$ for $t<0$, \eqref{wf1} can be rewritten as
$$
u^g(x,t)=\left \{
\begin{array}{cc}
f(t-x)+\int_{s=x}^tk(x,s)f(t-s)ds,&  x<t,\\
0,& x\geq t, \end{array} \right . $$
and similarly for \eqref{wf1d}.

\end{subsection}

\begin{subsection}{Solution representation on an interval}

Consider the IBVP on the interval $(0,l )$: 
\begin{eqnarray}
{u}_{tt}-{u}_{xx}+q(x){u}& = & 0,\ 0<x<l,\ t\in (0,T),\label{IBVP2a} \\
{u}(x,0)={u}_t(x,0) &= & 0,  \ 0<x<l,\nonumber \\
{u}(0,t) & = & h(t), t\in (0,T), \nonumber  \\
{u}(\ell ,t) &=  & 0, \ t\in (0,T).\label{IBVP2}
\end{eqnarray}
The homogeneous boundary condition at $x=\ell$ allows of a very convenient Goursat type representation of the solution $u$. Although this observation is perhaps well known, we first observed it in \cite{AZ}.

We extend $q$ to $(0,\infty)$ as follows: first evenly with respect to $x=l$, and then periodically. Thus $q(2n\ell \pm x)=q(x)$ for all positive
integers $n$.
Then the solution to \eqref{IBVP2a}-\eqref{IBVP2} on $(0,l)$ can be written as
$${u}^h(x,t)=\sum_{n\geq 0:\ 0\leq 2nl+x\leq t}
\left ( h(t-2nl-x)+\int_{2nl+x}^tk(2nl+x,s)h(t-s)ds\right )$$
\beq \label{fold}
-\sum_{n\geq 1:\ 0\leq 2nl-x\leq t}
\left (h(t-2nl+x)+\int_{2nl-x}^tk(2nl-x,s)h(t-s)ds \right ).
\eeq
This formula will be used in the proofs of { the shape and velocity controllability} in the  Section \ref{s3}. 

We now consider two IBVPs on the interval $(0,\ell )$  with Dirichlet and Neumann boundary conditions. First,
\begin{eqnarray}
{u}_{tt}-{u}_{xx}+q(x){u}& = & 0,\ 0<x<l,\ t\in (0,T),\label{IBVP2a} \\
{u}(x,0)={u}_t(x,0) &= & 0,  \ 0<x<l, \\
{u}(0,t) & = & h(t), \ t\in (0,T), \\
{u}_x(l ,t) &=  & 0, \ t\in (0,T). \label{IBVP2}
\end{eqnarray}
In this case, the solution to \eqref{IBVP2a}-\eqref{IBVP2} on $(0,l)$ can be written as
\begin{eqnarray}
{u}^h(x,t)& =& \sum_{n\geq 0:\ 0\leq 2nl+x\leq t}\left ( h(t-2nl-x)+\int_{2nl+x}^tk(2nl+x,s)h(t-s)ds\right )\nonumber \\
& &
+\sum_{n\geq 1:\ 0\leq 2nl-x\leq t}\left (h(t-2nl+x)+\int_{2nl-x}^tk(2nl-x,s)h(t-s)ds \right ).\label{fold2}
\end{eqnarray}

The following point-control result will be used in the Section 3. 
\begin{lemma}\label{pw}
 Let $\ep >0$, and let $T=l+\ep$.
Given $\phi \in H^1(0,l)$, there exists 
$h\in H^1(0,T )$, with $h(0)=0$ and with support  $[0,2\ep]$, such that the solution $\tilde{u}$
to \rq{IBVP2a}-\rq{IBVP2} satisfies
$\tilde{u}(l,T )=\phi (l)$. Furthermore, there exists a constant $C$ depending only on $l, q,\ep $ such that
$$\|h\|_{H^1(0,T)} \leq C||\phi ||_{H^1(0,l)}.$$
\end{lemma}
Proof: 
For any $\ep>0$, we define 
$$
L_{\ep}(t)=\left \{
\begin{array}{cc}
\ep^{-1/3}-\ep^{-4/3}|t|& |t|\leq \ep,\\
0 & |t|\geq \ep. 
\end{array}\right .
$$
Then as $\ep\to 0^+$, $L_{\ep}(0)\to \infty$ but $\| L_{\ep}\|_{L^2(\bR)}\to 0$. 
Let $\tilde{L}_\ep(t)=L_\ep(t-\ep).$
Then by \rq{fold2}, we have
$$
u^{\tilde{L}_\ep}(l,l+\ep)=2L_\ep(0)+\mbox{integral terms}.
$$
As $\ep \to 0^+$, the integral terms converge to zero, so there exists $\tilde{\ep}$ and $\al >0$ such that 
$u^{\tilde{L}_{\tilde{\ep}}}(l,l+\tilde{\ep})=\al$.
 We to set $h=\phi (l)\tilde{L}_{\tilde{\ep}}/\al$. Then we have 
$$\|h\|_{H^1(0,T)} \leq C ||{\phi }||_{H^1(0,l)}.$$
$\Box$

\ 

Finally, we discuss a solution representation when a Neumann control is applied at $x=l$. 
\begin{eqnarray}
{u}_{tt}-{u}_{xx}+q(x){u}& = & 0,\ 0<x<l,\ t\in (0,T),\label{IBVP9a} \\
{u}(x,0)={u}_t(x,0) &= & 0,  \ 0<x<l, \\
{u}(0,t) & = & 0, \ t\in (0,T), \\
{u}_x(l ,t) &=  & p(t), \ t\in (0,T).\label{IBVP9b}
\end{eqnarray}
Let $P(t)=-\int_0^tp(s)ds$.

Set $\tilde{q}(x)=q(\ell-x)$, and extend $\tilde{q}$ to $[0,\infty )$ by  $\tilde{q}(2k\ell \pm x)=\tilde{q}(x)$.
Define $w$  to be the solution to the Goursat-type problem
\begin{equation}\label{Gou2}
\left \{
\begin{array}{cc}
 w_{ss}(x,s)- w_{xx}(x,s)+\tilde{q}(x)w(x,s)= 0,& \ 0< x<s,\nonumber \\
w_x(0 ,s)  =  0, \
w(x ,x)  =  -\frac{1}{2}\int_0^x\tilde{q}_j(\eta )d\eta ,\ x<\ell .
\end{array} \right .
\end{equation}
In this case, the solution to the system \rq{IBVP9a}-\rq{IBVP9b} is given by 
\begin{eqnarray}
{u}^p(x,t)& =& \ P(t-\ell+x) +\int_{\ell-x}^t w( \ell-x,s) P(t-s) \, ds \nonumber \\
&&-P(t-\ell-x) - \int_{\ell+x}^t w( \ell+x,s) P(t-s) \, ds \nonumber\\
&&-P(t-3\ell+x) - \int_{3\ell-x}^t w( 3\ell-x,s) P(t-s) \, ds \nonumber\\
&&+P(t-3\ell-x) + \int_{3\ell+x}^t w( 3\ell+x,s) P(t-s) \, ds \nonumber\\
&&\dots \label{fold2bk}  \label{ndp}
\end{eqnarray}

\end{subsection}
\end{section}

\begin{section}{Proofs of shape and velocity control}\label{s3}
\begin{subsection}{Proof of Theorem \ref{thm3}}

We recall Problem 1:
\beq \lab{eqq}
u_{tt} - u_{xx} + q(x) u =0
\eeq
\beq \lab{bcc} 
\pa _x u_1(l,t)=f_1(t)
 \eeq
 \beq \lab{inc}
u_2(0,t)-u_1(0,t)=f_2(t), \ u_3(0,t)-u_1(0,t)=0, \ \sum_{j=1}^3 \partial u_j(0,t)=0,  
\eeq
with zero initial conditions and 
with the artificial vertex in the 
 the middle of the ring, where we  put the homogeneous Kirchoff--Neumann conditions: 
\beq \lab{ia}
u_2(a,t)-u_3(a,t)=0, \ \sum_{j=2}^3 \partial u_j(a,t)=0.  
\eeq

Recall  $T_*:=a+l.$ We now prove a slightly more precise version of 
Theorem \ref{thm3}.

{\bf Theorem 3$^*$.}

A) Let $T= T_*$. For any $\phi \in {\mH}^1_0,$ there exist $f_1 \in L^2(0,T)$ and $f_2 \in H^1_0(0,T)$ such that the solution $u$ to System \rq{eqq}-\rq{inc} satisfies $u(x,T)=\phi(x), x \in \Omega.$

B) Let $T=T_*$. For any $\psi \in {\mH},$ there exist $f_1 \in L^2(0,T)$ and $f_2 \in H^1_*(0,T)$ such that the solution $u$ to System \rq{eqq}-\rq{inc} satisfies $u_t(x,T)=\psi(x), x \in \Omega.$

Proof:

We start with the shape  controllability. We assume for now that $\phi_j(a)=0$, indicating the necessary modifications for the general case at the end. We set $f_2(t)=0$ for  $t\in (0,l)$.
We also put $y_j(t):=u_j(0,t+l), \, j=1,2.$ From \rq{inc} it follows that 
\beq \lab{g12}
y_1(t)=y_3(t), \quad f_2(t)=y_2(t-l)-y_1(t-l).
\eeq

First we will prove that the edges $e_2$ and $e_3$ are shape controllable with the help of
$f_1(t), t \in (0,a)$ and $f_2(t), t \in (l,T),$ or equivalently, with the help of
$y_1(t)$ and $y_2(t), t \in (0,a).$ Then we will demonstrate the edge $e_1$ is shape controllable with the help of
$f_1(t), t \in (a,T).$

By \rq{wf1d}, 
the functions $u_2$ and $u_3$ are trivial for $t<l$, and for $t>l$ they are presented by the formulas
\beq \lab{u23}
u_j(x,t)=y_j(t-l-x)+\int_x^{t-l} k_j(x,s) y_j(t-l-s) ds, \ j=2,3, \ t \in (l,T).
\eeq
In \rq{u23}, we now set $t=T$ and $u_j(x,T)=\phi_j(x),$  obtaining  Volterra equations of the second kind (which we will now refer to as a VESK)  for $y_j(t),\;t \in (l,T),\;j=2,3.$ 
Thus we have solved for $y_j\in H^1(0,a)$. The function $f_2(t)$ for $t\in (l,T)$ is now found from \rq{g12}. Since
$\phi_2(a)=\phi_3(a)=0$, we deduce from that $y_2(0)=y_3(0)=0,$ hence by \rq{g12}
we have $f_2(l)=0$. 
Furthermore, for $j=1, 2,3$ we have $\phi_j(0)=\phi_k(0)$, hence $y_j(a)=y_k(a)$, and hence $f_2(T)=0$. Thus $f_2\in H_0^1(0,T).$
  
It remains to solve for $f_1$, which requires the formula for $u_1$.
The function $u_1$ is described by a more complicated formula which depends whether or not  $l \geq a.$  

We assume first that $l \geq a.$ Then for $t \in (0,T)$ we have using \rq{fold2},\rq{ndp}
$$
 u_1(x,t)= -F_1(t-l+x)-\int_{l-x}^{t} w(l-x,s)F_1(t-s)ds $$
 $$ + F_1(t-l-x) +
 \int_{l+x}^{t}w(l+x,s)F_1(t-s)ds $$
 \beq \lab{u1}
 + y_1(t-l-x)+\int_x^{t-l} k_1(x,s) y_1(t-l-s) ds,
\eeq
where $F_1(t)=\int_0^tf_1(s)\,ds.$
  
  We substitute the solutions  in \rq{u23}, \rq{u1}
  to compute $ \pa u_j(0,t).$ 
  Then we apply the last equation of \rq{inc} to get 
   the following equation where $t \in (l,T):$
 $$ 0=-2f_1(t-l)+2\int_l^t\pa _x w(l,s)F_1(t-s)ds -2w(l,l) F_1(t-l) $$
 $$-2y'_1(t-l) - y'_2(t-l) +\int_0^{t-l}\pa _x k_1(0,s)y_1(t-l-s)ds     $$
  \beq \lab{f1}
 +\int_0^{t-l}\pa _x k_2(0,s)y_2(t-l-s)ds  +\int_0^{t-l}\pa _x k_3(0,s)y_1(t-l-s)ds. 
\eeq
This equation can be rewritten as a VESK for $f_1(t),\;t \in (0,a);$ thus we solve for $f_1$ on this time interval.

So, given the shapes $\phi_j(x)=u_j(x,T),\;j=2,3$ we have found the corresponding controls
$f_1(t),\;t \in (0,a)$ and $f_2(t),\;t \in (l,T).$ Now we put $t=T$ and $u_1(x,T)=\phi_1(x)$ into the equation \rq{u1}:
$$
\phi_1(x)= -F_1(a+x)-\int_{l-x}^{T} w(l-x,s)F_1(T-s)ds $$$$ + F_1(a-x) +
\int_{l+x}^{T}w(l+x,s)F_1(T-s)ds $$
\beq \lab{F1}
+ y_1(a-x)+\int_x^{a} k_1(x,s) y_1(a-s) ds, \ \ x \in (0,l).
\eeq
Note that $F_1(t)=0=y_1(t)$  for $t<0.$ Taking into account that $F_1(t)$ is known for $t \in (0,a),$ the only unknown terms in \rq{F1} are the first two on the right hand side. 
Thus we can solve this VESK for
$F_1(t), $ and thus
$f_1(t),$ for  $t \in (a,T).$

\vskip2mm

If $l < a,$ Equation \rq{u1} contains more terms, and hence so do the equations corresponding to \rq{f1}, \rq{F1}, and we will have to solve them in steps of size $l.$ Using \rq{fold2},\rq{ndp},  Equation \rq{u1} take the form:
$$
u_1(x,t)= -F_1(t-l+x)-\int_{l-x}^{t} w(l-x,s)F_1(t-s)ds $$
$$ + F_1(t-l-x) +\int_{l+x}^{t}w(l+x,s)F_1(t-s)ds $$
$$+F_1(t-3l+x)+\int_{3l-x}^{t} w(3l-x,s)F_1(t-s)ds $$
$$ - F_1(t-3l-x)- \int_{3l+x}^{t}w(3l+x,s)F_1(t-s)ds + \ldots $$
$$+ y_1(t-l-x)+\int_x^{t-l} k_1(x,s) y_1(t-l-s) ds$$
$$+ y_1(t-3l+x)+\int_{2l-x}^{t-l} k_1(2l-x,s) y_1(t-l-s) ds$$
\beq \lab{u11}
- y_1(t-3l-x)-\int_{2l+x}^{t-l} k_1(2l+x,s) y_1(t-l-s) ds + \ldots .
\eeq

Instead of equation \rq{f1} we get, for $t\in (0,T)$
 $$0= -2f_1(t-l)+2\int_l^t\pa _x w(l,s)F_1(t-s)ds -2w(l,l) F_1(t-l) $$
 $$-2y'_1(t-l) - y'_2(t-l) +\int_0^{t-l}\pa _x k_1(0,s)y_1(t-l-s)ds     $$ 
$$ +\int_0^{t-l}\pa _x k_2(0,s)y_2(t-l-s)ds  +\int_0^{t-l}\pa _x k_3(0,s)y_1(t-l-s)ds$$
$$+2f_1(t-3l) -2\int_l^t\pa _x w(3l,s)F_1(t-s)ds +2w(3l,3l) F_1(t-3l) $$ 
 \beq \lab{f11}
 +2y'_1(t-3l)-2\int_{2l}^{t-l}\pa _x k_1(2l,s)y_1(t-l-s)ds +2k_1(2l,2l)y_1(t-3l) + \ldots  .
 \eeq
 Here all the terms involving $y_j$ are known. Assume first that $t<3l$, so that the terms involving $f(t-(2n+1)l)$ for $n\geq 1$ all vanish. Then \rq{f11} becomes a VESK for $f_1(t-l)$, allowing us to solve for $f_1(t)$ for $t<2l$.
 Now assume $t\in (3l,5l)$. Then in \rq{f11},  the terms involving $f_1(t-3l)$ are now known, while  terms involving $f_1(t-(n+1)l)$ for $n\geq 2$ all vanish. This allows us to 
 solve for $f_1(t)$ for $t<4l$.  Iterating this argument, we solve for $f_1(t)$, $t\leq a.$

To find  $f_1(t),\;t \in (a,T)$,
we return to equation \rq{u11} and set $t=T, u_1(x,T)=\phi_1(x).$ We get
$$
\phi_1(x)= -F_1(a+x)-\int_{l-x}^{T} w(l-x,s)F_1(T-s)ds $$$$ + F_1(a-x) +
\int_{l+x}^{T}w(l+x,s)F_1(T-s)ds $$
$$+F_1(T-3l+x)+\int_{3l-x}^{T} w(3l-x,s)F_1(T-s)ds $$$$ - F_1(T-3l-x) -
\int_{3l+x}^{T}w(3l+x,s)F_1(t-s)ds + \ldots $$
$$+ y_1(a-x)+\int_x^{a} k_1(x,s) y_1(a-s) ds$$
$$+ y_1(T-3l+x)+\int_{2l-x}^{a} k_1(2l-x,s) y_1(a-s) ds$$
\beq \lab{F11}
- y_1(T-3l-x)-\int_{2l+x}^{a} k_1(2l+x,s) y_1(a-s) ds + \ldots .
\eeq
In this equation, $F_1(t)$ is known for $t<a$, so all terms are known except the first two terms on the right hand side. Solving for this VESK for $F_1(a+x)$ with $x\in (0,T-a),$ we find $f_1(t)$ for $t\in (a,T). $

For the case $\phi(a)\neq 0$, let $\ep\in (0, \min (l,a))$. 
Using an adaptation of Lemma \ref{pw}, we can find  a small ``bump function", $\tilde{f}_2\in H^1[0,T]$, supported in $(l-\ep ,l+\ep)$, such that the solution $v(x,t)$ to the system \rq{eqq}-\rq{inc} with $f_1=0, f_2=\tilde{f}_2$ satisfies
$v(a,T)=\phi(a)$, and $\| \tilde{f}_2\|_{H^1(0,T)}\leq C\| \phi\|_{\mH^1}$
for some $C$ that depends only on $\ep, q$.

The proof of Part A of Theorem 3$^*$ is completed as follows.
We  apply  the previous argument to attain the shape $\phi(x)-v(x,T)$ from controls $f_1\in L^2(0,T)$ and $f_2\in H^1_0(l,T)$. The controls
for Part A of the theorem will then be $f_1$ and $f_2+\tilde{f}_2$.

We now sketch the proof of Part B.  By \rq{u23}, one derives the following equations
 involving the target values of $(u_j)_t(x,T)$, $j=2,3$:
 \beq \lab{ghp}
 \psi_j(x)=y_j'(a-x)+\int_x^a k_j(x,s)y_j'(a-s)ds, \ x \in [0,a]. 
 \eeq
 Solving these VESK, we find 
 $y_j', f_2' \in L^2(l,T)$ and then find $y_j, f_2 \in H^1(l,T)$ from the formulae:
 \beq{} y_j(t)=\int_l^t y_j'(s)ds,\ j=1,2,3, \ f_2(t)=\int_l^t f_2'(s)ds, \ t \in [l,T].\label{gh2}
 \eeq
 Then we set $f_2(t)=0$ for $t<l$.
Mimicking the argument in part A, we find $f_1$ separately in cases $l\geq a$ and $l<a$. For $l\geq a$, we use \rq{f1} to solve for $f_1(t)$ for $t\in (a,T)$; then for $t< a,$ we differentiate
\rq{u1} in $t$, and then solve the resulting VESK with $(u_1)_t(x,T)=\psi_1(x).$ The case $l<a$ is similar.
Theorem 3$^*$ is proved.$\Box$

Note that the function $f_2(t)$ found in \rq{gh2}  does not necessarily satisfy $f_2(T)=0.$

\end{subsection}

\begin{subsection}
{Problem 2}

 Now we consider the same graph with three controls at the interior vertex:
\begin{eqnarray}
u_{tt} - u_{xx} + q(x) u &=& 0,\lab{eqq2}\\
\pa  u_1(l,t)&=&0,\lab{bcc} \\
u_3(0,t)-u_1(0,t)&=&f_3(t), \label{1c}\\
u_2(0,t)-u_1(0,t)& =& f_2(t), \label{2c}\\ 
\sum_{j=1}^3 \partial u_j(0,t)&=& f_1(t),  \lab{3c}\\
u(x,0)&=&u_t(x,0)=0.\label{ic}
\end{eqnarray}
We let 
$$T^*=\max (a,l ), \mbox{ and } T>T^*.$$
Then $T^*$ is the minimal time necessary for a wave generated at $x=0$ to reach both $x=l$ and $x=a$.
We begin with the following technical lemma.

\begin{lemma}\label{tl1}
Let $\ep \in (0,\min (a,l)).$
Given $\phi\in \mH^1_0$, there exist $f_1\in L^2(0,2\ep)$ and $f_2,f_3\in H^1_0(0,2\ep)$ such that the solution  $\tilde{u}(x,t)$ to the system \rq{eqq2}-\rq{ic} 
satisfies
$$\tilde{u}_1(l,T^*+\ep)=\phi_1(l) \mbox{ and }\tilde{u}_j(a,T^*+\ep)=\phi_j (a),\ j=2,3.$$
\end{lemma}
Proof:
Define $y_j(t)=u_j(0,t)$ for $j=1,2,3.$
We assume $l\geq a$, leaving the other case to the reader.
First, by Lemma \ref{pw}, there exists
 $y_1\in C_0^2(0,2\ep)$ that will generate a solution $u$ such that $u_1(l,T^*+\ep)=\phi_1(l)$.  We  
then extend $y_1$ trivially to $(0,T^*+\ep)$. We set $y_3=0.$
Next,  
by the proof of Lemma \ref{pw}, we can choose $y_2\in C_0^2(l-a, l-a+2\ep)$ such that  we have the associated wave 
${u}$ satisfying 
${u}_j(a,T^*+\ep)=\phi_j(a),$ $j=2,3.$ 
We then use \rq{1c},\rq{2c},\rq{3c} to determine functions ${f}_j$
for $j=1,2,3.$ $\Box$

\begin{thm}  

\ 

A) Let $T>T^*$. For any $\phi \in {\mH}_0^1,$ there exist $f_1 \in L^2(0,T)$ and $f_2,f_3 \in H^1_0(0,T )$ such that $u(x,T )=\phi(x), x \in \Omega.$

B) For any $\psi \in {\mH},$ there exist $f_1 \in L^2(0,T^*)$ and $f_2,f_3 \in H^1(0,T^*)$, with $f_2(0)=f_3(0)=0$ such that $u_t(x,T^*)=\psi(x), x \in \Omega.$
\end{thm}
Proof.
We begin by proving shape controllability. We will assume 
$l\geq a$, leaving the simple adaptation for $l<a$ to the reader. 
We consider first special case where $\phi (l)=\phi(a)=0.$ In this case, we can choose our controls so
that the desired shape can be attained in the optimal time 
 $T= l ,$ and there will be no reflection of waves.

Recall $u_j(0,t)=y_j(t),\, j=1,2,3;$  for $j=2,3$ we set $y_j(t)=0$ for $t<T-a$. We can write the solution $u$ in the form, for $t\in (0,T)$,
\beq \lab{uj}
u_j(x,t)=y_j(t-x)+\int_x^t k_j(x,s) y_j(t-s) ds, \ j=1,2,3, 
\eeq
so by \rq{g2m}
$$ \partial u_j(0,t)=-y'_j(t)+\int_0^tr_j(s)y_j(t-s)ds, \ r_j(s):=\pa k_j(0,s). $$
Here, for $j=1$, \rq{uj} holds for $x\in (0,l)$, while $j=2,3$
we have $x\in (0,a)$.
Note that these representations of the solution remain valid because there are no wave reflections, and this is guaranteed by our setting $y_2(t)=y_3(t)=0$ for $t<T-a$. 
Substituting into \rq{1c}-\rq{3c} we obtain the equations
\beq \lab{gf1}
y_3(t)-y_1(t)=f_3(t), \  y_2(t)-y_1(t)=f_2(t), 
\eeq
\beq \lab{gfj}
\sum_{j=1}^3 [-y'_j(t)+ \int_0^tr_j(t)y_j(t-s)ds]=f_1(t). 
\eeq
We put $t=T$ and $u_j(x,T)=\phi_j(x)$ in \rq{uj} to obtain  for $y_j(t),\;t \in (0,T),\;j=1,2,3:$ 
\beq \lab{ujf1}
\phi_j(x)=y_j(T-x)+\int_x^T k_j(x,s) y_j(T-s) ds, \ j=1,2,3, 
\eeq
$$  \mbox{ with }x \in [0,l] \mbox{ on }e_1, \; x \in [0,a] \mbox{ on } e_2,e_3. $$
These are Volterra integral equations of the second kind, so
we now solve for $y_j$.  For $j=1$, $\phi (l)=0$ implies 
$y_1(0)=0$. For $j=2,3$, our assumption
$y_j(T-s)=0$ for $s\geq a$ implies that $\phi_j(a)=0$, as desired, and $y_j(0)=0$. Thus by \rq{gf1}, $f_2(0)=f_3(0)=0.$ Furthermore, by continuity we have $\phi_1(0)=\phi_2(0)=\phi_3(0)$; since $k_j(0,s)=0$ for $j=1,2,3,$ we deduce $y_1(l)=y_2(l)=y_3(l),$
which in turn implies $f_2(l)=f_3(l)=0$, so $f_2,f_3\in H^1_0(0,T)$. By \rq{gfj}, we have $f_1\in L^2(0,T).$
This proves the shape controllability  of our graph in the special case. 

We now prove the shape controllability in the general case.
Let $\ep =T-T^*$, and we assume without loss of generality that $\ep$ satisfies the hypothesis in Lemma \ref{tl1}.
By translation in time, 
it suffices to construct $f_j$ supported in the time interval $(-\ep,T)$, given initial conditions
\beq u|_{t=-\ep}=u_t|_{t=-\ep}=0.\label{-ep}
\eeq
Using Lemma \ref{tl1},
there exist $\tilde{f}_1\in L^2(-\ep ,T)$, and
$\tilde{f}_j\in H_0^1(-\ep,T )$ for $j=2,3$, such that the solution $\tilde{u}(x,t)$ to the System \rq{eqq2}-\rq{3c} with initial conditions \rq{-ep} will satisfy 
$$\tilde{u}(l,T )=\phi (l),\ 
\tilde{u}(a,T )=\phi(a).$$
Then the function $\tilde{\phi}(x):=\phi (x)-\tilde{u}(x,T )$
is in $\mH^1_0$, and vanishes at $x=l$ and $x=a$.
Hence, by the special case proven above,  there exist $f_1\in L^2(-\ep,T )$ and $f_2,f_3\in H^1_0(-\ep,T )$ such that
the associated solution  satisfies $u(x,T)=\tilde{\phi}(x)$. Finally, the desired
controls to solve the general shape control problem will be
$f_j+\tilde{f}_j$, $j=1,2,3.$

\

To prove the velocity controllability
we proceed the same way. Instead of equations \rq{ujf1} we obtain
\beq \lab{ujf}
\psi_j(x)=y'_j(T-x)+\int_x^T k_j(x,s) y'_j(T-s) ds, \ j=1,2,3,
\eeq 
where for $j=1$ we have $\psi_1 \in L^2(0,l)$, while for 
$j=2,3$
 we have 
  $x\in [0,a], $
with $\psi_{j} \in L^2(0,a).$ Solving these equations, we find $y'_j,$ and then $y_j(t)=\int_0^t y_j'(s)\,ds,$ and then, using \rq{gf1}, \rq{gfj}, we find $f_j,\;j=1,2,3.$ Note that in this case, for $j=2,3$, we have $f_j \in H^1(0,T)$, and  $f_{j}(0)=0,$ but not necessarily $f_{j}(T)=0.$ $\Box$

\end{subsection}

\

We conclude this section with the following result on the well-posedness of Problems 1 and 2. The proof is based on the solution representations found in Section 2 together with the arguments found in this section, and is left to the reader.
\begin{prop}
Let $T>0$. Assume $q\in L^1(\Omega)$. 
Denote $H_*^1(0,T)=\{ f\in H^1(0,T):\ f(0)=0\}$.
Then,

A) for Problem 1, the mapping $(f_1,f_2)\mapsto (u,u_t)$ is a continuous map from $L^2(0,T)\times H_*^1(0,T)$ to $\mH_1^1\times \mH. $

B) for Problem 2, the mapping $(f_1,f_2,f_3)\mapsto (u,u_t)$ is a continuous map from $H^1_*(0,T)\times H^1_*(0,T)\times L^2(0,T)$ to $\mH^1\times \mH. $

\end{prop}
\end{section}
\begin{section}{Exact controllability via moment problems}\label{s4}

\begin{subsection}{Spectral preliminaries}

We now develop a spectral representation of the solution $u$ to prove { the controllability of our systems.}
Let $\{ (\omega^2_n,\f_n): n\geq 1\}$ be the eigenvalues and unit eigenfunctions for the spectral problem associated to
Problems 1, 2. Denote the restriction of $\f_n$ to edge $e_j$ by $(\f_n)_j$.

Thus
\beq \lab{sp}
 - \f_n'' + q(x) \f_n = \om_n^2 \f_n, \ x\in \cup_{j=1}^3 \, e_j \,, 
\eeq
\beq \lab{nc} 
(\f_n)_1' (l)=0,
\eeq
\beq \lab{kn}
(\f_n)_1(0)=(\f_n)_2(0)=(\f_n)_3(0), \ \sum_{j=1}^3  (\f_n)'_j(0)=0.   
\eeq
\begin{lemma}\label{specasy}
There exist positive constants $C_1,C_2,C_3$ such that 
the eigenvalues $\{ \omega_n^2\}_1^{\infty}$ associated \rq{sp}-\rq{kn}
satisfy 
$$C_1 n \leq |\omega_n| +1 =C_2n,$$
and the associated unit norm eigenfunctions satisfy
$$|\f_n (l)|\leq C_3,\ |(\f_n)_j'(0)|\leq C_3 n.$$
\end{lemma} 
 For more information about spectral problems on metric graphs
see, e.g. \cite{BK}, \cite{Ku2} and \cite{AN}. 

\end{subsection}
\begin{subsection}{Proof of Theorem \ref{thm1}}
We first briefly outline the proof. After developing a spectral solution to the IBVP in Problem 1, we use shape and velocity controllability to prove solvability of a pair of moment problems. Then by a simple argument involving the symmetry of the sine and cosine functions, we solve the moment problem associated to exact controllability.

For  just this subsection, for readability, we will make the following notational change: the boundary control with be denoted $f$ instead of $f_1$, and the interior control  will be denoted $h$ instead of $f_2$.

Let $\phi \in {\mH}^1_0$ and  $\psi \in {\mH}.$ Let $T=T_*.$
By Theorem \ref{thm3} we have shape and velocity controllability.

The solution of the IBVP \rq{eqq}--\rq{inc}, \rq{ia} can be presented in a form of series
\beq \lab{ser}
u(x,t)=\sum_{n=1}^{\infty} a_n(t)\f_n(x).
\eeq
To find the coefficients $a_n$ we integrate by parts the identity
$$0=\int_0^T \int_{\Omega} [u_{tt}-u_{xx}+q(x)u]\,\f_n(x)\mu(t)\,dxdt$$
with arbitrary $\mu \in C^2_0[0,T].$ We obtain
the initial value problem
$$a_n''(t)+\omega_n^2a_n(t)=f(t)(\f_{n})_1(l)
+h(t)(\f_n)'_2(0), \ a_n(0)=a_n'(0)=0.
$$
By variation of parameters, we get
\beq \lab{an}
a_n(T)=\int_0^T [f(t)(\f_{n})_1(l)+h(t)(\f_n)'_{2}(0)]\,\frac{\sin \om_n(T-t)}{\om_n}\,dt.
\eeq
We assume here that $\omega_n >0.$ If  $\omega_n =0,$ instead of $\frac{\sin \omega_n(T-t)}{\om_n}$
we put $T-t,$ and if  $\om_n <0,$ we put $\frac{\sinh \om_n(T-t)}{\om_n}.$

First, we use the shape controllability to derive a moment equation. 
Integrating by parts the term with $h$ and using that $h \in H^1_0(0,T),$ we rewrite
\rq{an} in the form 
\beq \lab{an1}
 a_n(T)\,\om_n=\int_0^T [f(t)\,\al_n \sin \om_n(T-t)-h'(t)\,\be_n\cos \om_n(T-t)]\,dt
 \eeq
with $\al_n= (\f_{n})_1(l)$ and  $\be_n = (\f_n)'_{2}(0)/ \om_n$.

Next, we derive a moment equation from velocity control. 
Now we differentiate \rq{an} with respect to $T$,  and then integrate by parts the term with $h$, using $h\in H^1(0,T)$ and $h(0)=0$. We obtain
\beq \lab{and}
{\dot a}_n(T)=\int_0^T [f(t)\,\al_n \cos \om_n(T-t)+h'(t)\,\be_n\sin \om_n(T-t)]\,dt.
\eeq
Using the classical relations between the Fourier method, the moment problem and control theory (see, e.g. \cite[Ch. III]{AI}), we can reformulate the statements Theorem \ref{thm3} as follows:

(a1) For any sequence $\{a_n\om_n\} \in \ell^2,$ there exist $f_0 \in L^2(0,T),\; h_0 \in
H^1_0(0,T)$ such that
\beq \lab{anm}
a_n\,\om_n=\int_0^T [f_0(t)\,\al_n \sin \om_n(T-t)-h_0'(t)\,\be_n\cos \om_n(T-t)]\,dt.
\eeq

(b1) For any sequence $\{b_n\} \in \ell^2,$ there exist $f_1 \in L^2(0,T),\; h_1 \in
H^1(0,T)$ with $h(0)=0$ such that
\beq \lab{bn}
b_n=\int_0^T [f_1(t)\,\al_n \cos \om_n(T-t)+h_1'(t)\,\be_n\sin \om_n(T-t)]\,dt.
\eeq
Now we are ready to complete the proof  of Theorem 1. We extend the functions $f_0,f_1$ and $h'_{0},h'_1$ from the interval $[0,T]$ to $[0,2T]$ such that
for $t \in [0,T]$ we have
$$f_0(t)=-f_0(2T-t), \ h'_1(t)=-h'_1(2T-t), \ f_1(t)=f_1(2T-t), \ h'_0(t)=h'_0(2T-t)$$
and put
$$ f(t)=\frac{f_0(t)+f_1(t)}{2}, \quad h'(t)=\frac{h'_0(t)+h'_1(t)}{2}.$$
One can check that
\beq \lab{anf}
a_n\,\om_n=\int_0^{2T} [f(t)\,\al_n \sin \om_n(T-t)-h'(t)\,\be_n\cos \om_n(T-t)]\,dt.
\eeq
\beq \lab{bnf}
b_n=\int_0^{2T} [f(t)\,\al_n \cos \om_n(T-t)+h'(t)\,\be_n\sin \om_n(T-t)]\,dt.
\eeq

{
Solvability of the moment problem \rq{anf}, \rq{bnf} is equivalent to solvability of
another  moment problem:
\beq \lab{aef}
c_n=\int_0^{2T} [f(t)\,\al_n -i h'(t)\,\be_n]\,e^{i \om_n(T-t)}\,dt.
\eeq
\beq \lab{bef}
d_n=\int_0^{2T} [f(t)\,\al_n +ih'(t)\,\be_n]\,e^{-i \om_n(T-t)}\,dt,
\eeq
for any sequences $\{c_n\},\, \{d_n\} \in \ell^2$ and the same spaces of functions $f,h' \in L^2(0,2T).$ Clearly, in these equalities one can change $T$ to $2T$ in the arguments of the exponentials. Finally, from $\exp\{\pm i \om_n(2T-t) \}$ we can switch back to $\sin \om_n(2T-t) ,\,  \cos \om_n(2T-t) $.
This allows us to claim solvability of the moment problem \beq \lab{an2}
a_n\,\om_n=\int_0^{2T} [f(t)\,\al_n \sin \om_n(2T-t)-h'(t)\,\be_n\cos \om_n(2T-t)]\,dt.
\eeq
\beq \lab{bn2}
b_n=\int_0^{2T} [f(t)\,\al_n \cos \om_n(2T-t)+h'(t)\,\be_n\sin \om_n(2T-t)]\,dt.
\eeq
}
We put $h(t)=\int_0^t h'(s)\,ds$ and check that
$$2h(2T)= \int_0^{2T} h_0'(s)\,ds + \int_0^{2T} h_1'(s)\,ds=2\int_0^{T} h_0'(s)\,ds=0,
$$ 
hence $h \in H^1_0(0,2T).$ Also, evidently, $f\in L^2(0,2T).$
Now integrating by parts, we can rewrite \rq{an2}, \rq{bn2} in the form
\beq \lab{af}
a_n\,\om_n=\int_0^{2T} [f(t)\,\al_n +h(t)\,\be_n \om_n ]\, \sin \om_n(2T-t)\,dt.
\eeq
\beq \lab{bf}
b_n=\int_0^{2T} [f(t)\,\al_n +h(t)\,\be_n \om_n] \, \cos \om_n(2T-t)\,dt.
\eeq
Thus for arbitrary sequences $\{a_n\om_n\},\; \{b_n\}  \in \ell^2,$ there exist $f \in L^2(0,2T),\; h \in H^1_0(0,2T),$ satisfying equalities \rq{af}, \rq{bf}. This is equivalent to the exact controllability of Problem 1. $\Box$

We remark in passing that by Lemma \ref{specasy}, there exists a constant $C>0$ such that 
$\al_n<C,  \be_n<C$, ensuring the appropriate regularity for \rq{ser}.

\end{subsection}
\begin{subsection}{Proof of Theorem \ref{thm2}}


Now we prove the exact controllability in time $2T,$
using the shape and velocity controllability in time $T$.
Using the same spectral representation of the function $u(x,t)$ as in Problem 1, we obtain 
\beq \lab{3an}
a_n(T)=\int_0^T [f_3(t)\,\f_n(0)-\sum_{j=1}^2 f_j(t)\,(\f_n)'_{j}(0)]\,\frac{\sin \om_n(T-t)}{\om_n}\,dt.
\eeq
In what follows, the controls solving the shape control problem will be denoted by 
$f_j=f_{j0}$, and the controls solving the velocity control problem will be denoted by $f_j=f_{j1},$ for $j=1,2,3$
We apply integration by parts for the terms with $f_1,f_2$ in \rq{3an}, using that $f_{10},f_{20} \in H^1_0(0,T).$ Setting $a_n=a_n(T)$, the shape controllability result can be formulated as solvability of the moment problem
\beq \lab{3a}
a_n \om_n=\int_0^T [f_{30}(t)\,\ka_{n3}\,\sin \om_n(T-t)+\sum_{j=1}^2 f'_{j0}(t)\,\ka_{nj}\,\cos \om_n(T-t)]\,dt,\ \forall n, 
\eeq
for arbitrary sequence $\{a_n \om_n\} \in \ell^2$ by the functions $f_{30},f_{10}',f_{20}' \in L^2(0,T).$ Here $$\ka_{n3}=(\f_n)_3(0),\ \ka_{nj}=\frac{(\f_n)'_{j}(0)}{\om_n}\, \ j=1,2.$$
 
Differentiating \rq{3an} with respect to $T$, and  then integrating by parts while using $f_{11}(0)=f_{21}(0)=0$, we reformulate the 
the velocity controllability result  as solvability of the moment problem
\beq \lab{3bn}
b_n =\int_0^T [f_{31}(t)\,\ka_{n3}\,\cos \om_n(T-t)-\sum_{j=1}^2 f'_{j1}(t)\,\ka_{nj}\,\sin \om_n(T-t)]\,dt, \ \forall n,
\eeq
for arbitrary sequences $\{b_n\} \in \ell^2$ by the functions $f_{11},f_{21}',f_{31}' \in L^2(0,T).$ 

Now we extend the functions $f_{30},f'_{11},f'_{21}$ in the odd way with respect to $t=T$
from the interval $[0,T]$ to $[0,2T], $ and the functions $f_{31},f'_{10},f'_{20}$ --- in the even way. Define the functions 
$$f_3=\frac{f_{30}+f_{31}}{2}, 
\ f'_{1}=\frac{f'_{10} + f'_{11}}{2}, \ 
f'_{2}=\frac{f'_{20}+f'_{21}}{2},	$$
and observe that these functions solve the moment problems
\beq \lab{3a1}
a_n \om_n=\int_0^{2T} [f_{3}(t)\,\ka_{n3}\,\sin \om_n(T-t)+\sum_{j=1}^2 f'_{j}(t)\,\ka_{nj}\,\cos \om_n(T-t)]\,dt
\eeq
\beq \lab{3b1}
b_n =\int_0^{2T} [f_{3}(t)\,\ka_{n3}\,\cos \om_n(T-t)-\sum_{j=1}^2 f'_{j}(t)\,\ka_{nj}\,\sin \om_n(T-t)]\,dt
\eeq
for arbitrary sequences $\{a_n \om_n\}, \; \{b_n\} \in \ell^2$. Also, defining $f_{j}(t)=\int_0^t f'_{j}(s)\,ds,$ for $j=1,2,$ it's not hard to prove $f_j\in H^1_0(0,T)$ for $j=1,2.$
Thus 
these moment problems can be rewritten in the form
\beq \lab{3a2}
a_n \om_n=\int_0^{2T} [f_{3}(t)\,\ka_{n3}\ -\sum_{j=1}^2 f_{j}(t)\,\ka_{nj}\omega_n]\,\sin \om_n(T-t)\,dt
\eeq
\beq \lab{3b2}
b_n =\int_0^{2T} [f_{3}(t)\,\ka_{n3}-\sum_{j=1}^2 f_{j}(t)\,\ka_{nj}\omega_n]
\,\cos \om_n(T-t)\,dt,
\eeq
Having solved these moment problems is equivalent to proving exact controllability in time $2T$. The 
proof of Theorem \ref{thm2} is complete.$\Box$

\end{subsection}

\end{section}

\begin{section}{Non-controllability for one control}\label{s5}

The main purpose of this section is to show that if the number of controls is reduced in Problem 1, then the resulting system might  not be controllable.

\begin{subsection}{Lack of approximate controllability for $q=0$.}
For this section, for Problem 1, we use $f_1$ for the boundary control, and $f_2$ for the interior control.

We now show that if $q=0$ and either of $f_1,f_2$ is set to zero, then  the system described in
Problem  1 is not  approximately controllable for any $T>0$.

Suppose first that $f_1=0.$ Then equality \rq{an} takes the form
\beq \lab{an0}
a_n(T)=\int_0^T f_2(t)(\f_n)'_{2}(0)\,\frac{\sin \om_n(T-t)}{\om_n}\,dt.
\eeq
For $q=0,$ there is an eigenfunction of our problem which is constant on
$\Om.$ We will denote it by $\f_1.$ From \rq{an0} it follows that $a_1^{f_2}(T)=0.$  Therefore $u^{(0,f_2)}$ $(\cdot,T)$ 
is orthogonal to $\f_1$ for all $f_2$ and any $T>0.$

Suppose next that $f_2=0.$ Then by symmetry,  $u_2^{(f_1,0)}(x,t)=u_3^{(f_1,0)}f(x,t), \, x \in (0,a), \, t>0,$  and the system is clearly not approximately controllable for any $T>0$.
	
	We provide also another argument for proving this fact which will be used in the next section.
If $f_2=0,$  equality \rq{an} takes the form
\beq \lab{anf1}
a_n(T)=\int_0^T f_1(t)\,(\f_{n})_1(l)\,\frac{\sin \om_n(T-t)}{\om_n}\,dt.
\eeq
There is an infinite number of the eigenfunctions vanishing on $e_1:$ for each $k=1,2,\ldots $,  there exists $n_k >1$ such that
$$
\varphi_{n_k} (x)=
\left \{ 
\begin{array}{cc}
0, & x\in e_1, \\
\sin (\pi k x/a), & x\in e_2\\
-\sin (\pi k x/a), & x\in e_3,
\end{array}
\right .
$$
Thus  $a_{n_k}^{f_1}(T)=0, \; \forall n_k.$ 
From \rq{anf1} it follows that $u^{(f_1,0)}(\cdot,T)$ is orthogonal to all $\f_{n_k}$ for all $f_1$ and any $T>0.$

\end{subsection}
\begin{subsection}{Lack of exact controllability for general $q$}
	
	 In this section we demonstrate another kind of negative controllability result.
		We prove that our system (related to Problem 1) is never exactly controllable (for arbitrary continuous potential) if we use only one controller that produces either boundary controls $f_1$ or interior controls $f_2$. If we use only $f_1$, the problem is more simple. As we demonstrated in the previous section, for trivial $q$ there are an infinite number of eigenfunctions supported on the ring. Then, for any $q$ and any given $n$, there is an eigenfunction $\f_{m_n}^q$ arbitrary close to one those eigenfunctions, so
		$|(\f_{m_n}^q)(l)| < 1/n.$ Therefore, the family $\{\f_{1,n}(l) \sin \om_n(T-t) \} $
			is not a Riesz sequence in $L^2(0,T)$ for any $T>0,$ and the moment problem
	\beq \lab{anft}
	a_n(T) \om_n =\int_0^T f_1(t)\,(\f_{n})_1(l)\,\sin \om_n(T-t)\,dt.
	\eeq
cannot be solvable for arbitrary sequence $\{a_n(T)\, \om_n\} \in \ell^2$ (see \cite[Sec. I.2]{AI}).

\

In the case when we use only the interior control $f_2$ we need a different argument to 	demonstrate that the corresponding 	moment problem is not solvable. We prove that the sequence of the eigenfrequencies is not uniformly discrete. It is convenient now to use a slightly different notations.

Consider a lasso with $e_1$ of length $l$ and ring, $e_2,$ of circumference $a$. For this section, we parametrize
$e_1$ with 
the boundary vertex at $x=0$. Then we parametrize  the ring by
$x\in (0,a)$ with $x=0$ and $x=a$ identified on the graph.

\begin{prop}\label{gap}
The sequence of the eigenfrequencies $\{\omega_n\}$ of the wave equation $u_{tt}-u_{xx}+q(x)u$
on this lasso graph is not separated:
$$\inf_{m\neq n}|\omega_m-\omega_n|=0.
$$
\end{prop}
Proof of proposition: We prove the result
for  $q=0$;  the case $q\neq 0$ will then follow  by a mini-max argument found in \cite{BKS}.
	
	Let $\{ (\omega^2_n,\f_n): n\geq 1\}$ be the eigenvalues and unit eigenfunctions for the 	
	following spectral problem 
	\beq \lab{sp2}
	- \f''  = \om^2 \f, \ x \in \cup_{j=1}^2 e_j\,,
	\eeq
	\beq \lab{nc2} 
	\f_1' (0)=0,
	\eeq
	\beq \lab{cont2}
	\f_1(l)=\f_2(0)=\f_2(a),
	\eeq
	\beq \lab{kn2}
	\ -\f'_1(l)+\f_2'(0)-\f_2'(a)=0.   
	\eeq
	
Because the quadratic form associated with this problem is non-negative, we have $\om^2\geq 0.$ We also note that $\om =0$ is a frequency with eigenspace spanned by  $\f =1$. Henceforth, we will assume $\om >0.$

Next, we consider the case where $\f_1=0$. Thus $\f_1(l)=\f'(l)=0$. It follows that $\f_2(x)=\sin (2n x \pi/a)$, and the frequency spectrum is 
$$\Lambda_1:= \{ 2n\pi /a,\ n\in \bN\}.$$ 
Next, we note that if $\f_2=0$, this \rq{cont2} and \rq{kn2} imply $\f =0$.

Next, we consider eigenfunctions non-trivial on $e_1$.
On $e_1$, by \rq{sp2} and \rq{nc2}, we can set
$$\f _1(x)=\cos (\om x).$$
Hence 
\beq 
\f_1(l) =\cos (\om l),\ \pa \f_1(l)=\om \sin (\om l).\label{e1}
\eeq
On the interval $(0,a)$,  we have 
$$
\f_2 (x)=A\cos (\om x)+B\sin (\om x).
$$
We apply $\f_1(l)=\f_2(0)$ to get
\beq
\f_2 (x)=\cos (\om l)\cos (\om x)+B\sin (\om x),\label{f2}
\eeq
so that 
\beq 
\f_2' (0)=B\om , \ \f_2'(a)=-\om \cos (\om l)\sin (\om a)+B\om \cos (\om a).
\label{kn0}
\eeq
Then \rq{kn2}  implies
\beq
0=\sin (\om l)+B+ \cos (\om l)\sin (\om a)-B \cos (\om a).\label{B}
\eeq
{\bf Case 1: $a/l\in \bQ.$}
Thus $a/l=p/q,$ with $p,q\in \bN.$ Consider 
$$
\Lambda_2:= \{ 2\pi mp/a: m\in \bN\}.
$$
Then  it is easily verified that $\om \in \Lambda_2$ satisfies \rq{B} and that \rq{cont2}, \rq{kn2} hold, so $\Lambda_2$ forms a frequency set with eigenfunctions linearly independent of those associated to $\Lambda_1$. Since $\Lambda_2\subset \Lambda_1$, we have an infinite set of frequencies of multiplicity two. Thus proves the proposition in this case.

{\bf Case 2: $a/l\notin \bQ.$}
Suppose first $1=\cos (\om a)$. Then \rq{B} implies 
$\sin (\om l)=0.$ Then we have both $\om =2n\pi /l$ and 
$\om=m\pi /a$ for some $ m, n  \in \mathbb{N}, $ implying $a/l\in \bQ$, a contradiction. 

Hence we have $1\neq \cos (\om a)$.
Thus, 
 we solve for $B$ in \rq{B}:
$$B=-\frac{ \sin (\om l)+ \cos (\om l)\sin (\om a)}
{1-\cos (\om a)}.$$
We now use $\f_2(a)=\cos (\om l)$ and \rq{f2} to get 
\beq 
0 =2 \cos (\om l)\cos (\om a)-2\cos (\om l)- \sin (\om l) \sin(\om a)\label{char}
\eeq
We consider first the case $a<l.$
We will use a Diophantine argument.
Let $$\la_n=\frac{2\pi n}{l}.$$
Denote $|||x|||$ the distance from $x$ to the lattice
$\{ 2\pi n:\ \in \bN\} $.

\begin{lemma} \label{dio}
Suppose $a/l$ is irrational. Then there exists $C>0$ and an infinite
subsequence of $\{ \la_n\}$ satisfying
$$|||\la_n a |||\leq C/n.$$
\end{lemma}

The proof of the lemma is deferred to the end of the section.
In what follows, we will assume $n$ belongs to the subset of $\bN$ specified in the  lemma.
Let
\beq
g(\om )=2\cos^2(\om l)-2\cos (\om l)-\sin^2 (\om l).\label{defg}
\eeq
Thus $g$ vanishes to second order on the set $\{ 2\pi n/l\}.$
Let $f(\om )$ be the right hand side of \rq{char}. We will use Rouche's Theorem to show that $f$ has two roots close to
$\{ 2\pi n/l\}.$

To this end,
 let $C$ be as in the statement of that lemma. Also, we will assume $n$ satisfies 
\beq
\frac{n}{C} > 100 \ln (n).\label{log}
\eeq
Fix any such $n$.

Suppose $\om$ is now on a circle of radius $\frac{1}{l\ln (n)}$ centered at $2\pi n/l$
in $\bC$. In particular, we set $\om =2\pi n/l+z/l$, with 
$|z|=1/\ln (n).$
Then by periodicity of sine and cosine along with their Taylor expansions, we have 
\begin{eqnarray}
g(2\pi n/l+z/l) & = & 2\cos^2(z)-2\cos (z)-\sin^2(z)\nonumber \\
& = & -2z^2+O(z^4).\label{g}
\end{eqnarray}
Next, we estimate $f-g$ on the same circle.
\begin{eqnarray}
f(\om )-g(\om ) & = &2\cos (\om l)(\cos (\om a)-\cos (\om l))
+\sin (\om l)(\sin (\om l)-\sin (\om a))\nonumber \\
& = & -4\cos (\om l)\sin (\om (l+a)/2)\sin (\om (a-l)/2) \nonumber \\
& & +2\sin (\om l)\sin (\om (l-a)/2)\cos (\om (l+a)/2) \nonumber \\
& = & 2\sin (\om (l-a)/2)\big (\sin (\om l) \cos (\om (l+a)/2) \nonumber \\
&& +2\cos (\om l)\sin (\om (l+a)/2)\big ). \label{comp}
\end{eqnarray}

We have for $\om =\la_n+z/l$ , and by Lemma \ref{dio}, there exists
 $m=m(n)\in \bZ$ such that
$\la_n(l-a)=2m\pi +O(1/n)$ and
\begin{eqnarray*}
\om (l-a)/2 & = & (\la_n +z/l)(l-a)/2\\
& = & (2m\pi +O(1/n))/2+\frac{z}{2}(1-a/l),\\
\om (l+a)/2 & = & \om ((a-l)/2+l)\\
& = &(\la_n +z/l) ((a-l)/2+l)\\
& = & (-2m\pi +O(1/n))/2 +2n\pi +\frac{z}{2}(1+a/l).
\end{eqnarray*}
The following estimates hold
\begin{eqnarray*}
\sin (\om (l-a)/2) & = &(-1)^m\frac{z}{2}(1-a/l)+ O(1/n)+O(z^3)\\
\cos (\om (l+a)/2))& = & (-1)^m(1+O(z^2)+O(1/n))\label{cos}\\
\sin (\om (l+a)/2))& = &(-1)^m \frac{z}{2}(a/l+1)+ O(1/n)+O(z^3)\\
\cos (\om l) & = & 1+O(z^2)\\
  \sin (\om l) & = & z+O(z^3).
\end{eqnarray*}
Inserting this estimate into \rq{comp} give 
\beq
f(\om )-g(\om )={z^2}(1-a/l)(2+a/l)+O(z^3)+O(1/n).\label{f-g}
\eeq

Finally, {since $0<a/l<1$},
$$
|(1-a/l)(2+a/l)|<2.
$$
{If $n$ is sufficiently large,}
then by \rq{log}, \rq{g} and \rq{f-g}, we have 
$$|f(\om )-g(\om )|< |g(\om )|,$$
from which we conclude by Rouch\'e's Theorem that $f$ has two zeros in the disk
of radius $1/(l\ln (n))$ centered at $2\pi n/l.$

We now consider the case $l<a$.  Let
\beq
h(\om )=2\cos^2(\om a)-2\cos (\om a)-\sin^2 (\om a).\label{defh}
\eeq
Thus $h$ vanishes to second order on the set $\{ 2\pi n/a\}.$
We can use Rouche's Theorem to show that $f$ has two roots close to
$\{ 2\pi n/a\};$
the proof is parallel to the case $a<l$ and is left to the reader.

Thus we have proven the proposition in the case where $a/l$ is irrational.  $\Box$

\begin{cor}
If  $f_1=0$ in
the system \rq{eqq1}-\rq{inc1}, then the system is not exactly  controllable for any  $T>0,$ i.e. the conclusions of Theorem \ref{thm1} will fail for any $T>0$.
\end{cor}
Proof: Let $T>0$. 
By the proposition and Lemma \ref{specasy}, standard arguments in \cite{AI} show that
$\{ \om_n^{-1}(\f_n)_2'(0)e^{i\om_n t}\}$ does not form a Riesz sequence on $L^2(0,T) $, from which the corollary follows \cite[Theorem III.3.10]{AI}.$\Box$

Proof of Lemma \ref{dio}: 
The inequality above is equivalent to the existence of $m\in \bZ$ such that 
$$\left|\frac{2\pi na}{l}-2\pi m\right|\leq C/n,$$
i.e. 
$$\left|\frac{a}{l}-\frac{ m}{n}\right|\leq \frac{(C/2\pi)}{n^2}.$$
But it is a well known from the theory of Diophantine equations that the last inequality is satisfied by infinitely
many integer pairs $(m,n)$ with $2\pi /C=\sqrt{5}.$ $\Box$

\

We have proved that the system \rq{eqq1}--\rq{inc1} is not exactly  controllable for any  $T>0$ under the action of only one control, either $f_1$ or $f_2.$ From the moment problems
\rq{af}, \rq{bf} it follows that in this case the system still may be spectrally controllable for $T > 2L, \; L:=l+2a,$ if the spectrum of the problem \rq{sp}--\rq{kn} is simple and either $(\f_n)_1(l) \neq 0$ or $(\f_n)_2'(0) \neq 0, \ \forall n.$

 We give here only a sketch of the proof. Under pretty mild conditions on the regularity of $q$ the uniform density of the eigenfrequences $\om_n$ is equal to $L/\pi$ (see \cite[Sec. 2.1]{AN}, \cite[Sec. 3.7]{BK}). Then Theorem 3(ii) from \cite{AM} implies that the family
$\{\sin \om_nt,\, \cos \om_nt\}$ is minimal in $L^2(0,T)$ for $T>2L,$ and the system \rq{eqq1}--\rq{inc1} is spectrally  controllable in such time intervals.

Certainly, this controllability is very unstable with respect to small perturbations of the system parameters $a,l,q.$ A stable controllability is guaranteed only by two controls as described in Theorem 1.

\end{subsection}
\end{section}

\vskip3mm

\noindent {\bf  Acknowledgments}\\
The research of Sergei Avdonin was  supported  in part by the National Science Foundation,
	grant DMS 1909869. 
The research of Yuanyuan Zhao was supported by the National Science Foundation Graduate Research Fellowship under Grant No. 1242789.

\end{document}